\documentclass[10pt]{siamltex}
\usepackage{amsmath, amssymb, latexsym, array}

\usepackage{srcltx}
\usepackage{graphicx}
\usepackage{pst-all}

\newtheorem{remark}{Remark}[section]
\newtheorem{example}{Example}[section]
\newcommand{\bd}{\partial}
\newcommand{\Grad}{\nabla}

\renewcommand{\div}{\operatorname{div}}

\newcommand{\norm}[1]{\|#1\|}

\def\pd#1#2{\frac{\partial#1}{\partial#2}}

\def\mdiv{\div\,}

\def\vcurl{\mathbf{curl\,}}

\def\vrot{\mathbf{rot\,}}

\def\rot{\mathrm{rot\,}}
\def\OM{\Omega}
\def\R{\mathbb{R}} \def\RR {\mathbb R}

\def\ba{{\bf a}}
 
\def\bg{{\bf g}}

\def\bn{{\bf n}}
\def\bq{{\bf q}}

\def\bu{{\bf u}}
\def\bv{{\bf v}}
\def\bx{{\bf x}}

\def\bH{{\bf H}}
\def\bS{{\bf S}}

\def\bV{{\bf V}}

\def\cE{\mathcal E} \def\cG{\mathcal{G}}
\def\cK{\mathcal K}
\def\cR{\mathcal{R}}

\def\cT{\mathcal T}

\def\Poly{\mathcal P}

\def\balpha{\mbox{\boldmath{$\alpha$}}}
\def\bphi{\mbox{\boldmath{$\phi$}}}
\def\bsigma{\mbox{\boldmath{$\sigma$}}}

\def\bXi{\mbox{\boldmath{$\Xi$}}}

\def\bPi{\mbox{\boldmath{$\Pi$}}}
\def\bPsi {\mbox{\boldmath{$\Psi$}}}

\begin{document}

\title{Mixed Virtual Volume  Methods for Elliptic Problems}
\author{
Gwanghyun Jo\thanks{Department of Mathematics, Kunsan National University, 558 Daehak-ro, Gunsan-si, Jeollabuk-do, Republic of Korea  54150. This author is supported by  National Research Foundation of Korea(NRF) grant, contract No. 2020R1C1C1A01005396.
{\tt
email:gwanghyun@kunsan.ac.kr }}
\and
 Do Y. Kwak\thanks{Department of mathematical Sciences, Korea Advanced Institute of Science and Technology, 291 Daehak-ro, Daejeon, Republic of Korea 34141. This author is supported by NRF grant, contract No. 2021R1A2C1003340.
  {\tt email:kdy@kaist.ac.kr}} }

 \maketitle

\begin{abstract}
We develop a class of mixed virtual volume methods for elliptic problems on  polygonal/polyhedral grids. Unlike the  mixed virtual element methods introduced in \cite{brezzi2014basic,da2016mixed}, our methods are reduced to symmetric, positive definite problems for the primary variable without using Lagrangian multipliers.
  We start from the usual way of changing the given equation into a mixed system using the Darcy's law, $\bu=-{\cal
 K} \nabla p$.  By integrating the system of equations with some
judiciously chosen  test spaces  on each element,  we define new mixed virtual volume methods of all orders. We
show that these new schemes are equivalent to the nonconforming virtual element methods for the primal variable $p$.

Once the primary variable is computed solving the symmetric, positive definite system, all the degrees of freedom
for the Darcy velocity are locally computed.
Also, the $L^2$-projection onto the polynomial space is easy to compute.
Hence our work opens an easy way to compute  Darcy velocity  on the polygonal/polyhedral grids.
For the lowest order case, we give a formula
to compute a Raviart-Thomas space like representation which satisfies the conservation law.

   An optimal error analysis is carried out and numerical results are presented  which support the theory.
\end{abstract}


\noindent {\bf AMS(MOS) subject classifications.} 65N15, 65N30.

\pagestyle{myheadings} \thispagestyle{plain} \markboth{
Gwanghyun Jo and Do Y. Kwak}{Mixed Virtual Volume Methods for Elliptic Problems}


\section{Introduction}

The virtual element method (VEM), introduced by Beir\~{a}o da Veiga, et al. \cite{beirao2013basic}, is a  generalization of the conventional finite element method to general polygonal (or polyhedral) meshes, where thorough error analysis and numerical tests for more general cases for elliptic problems were developed in \cite{beirao2013basic,  ahmad2013equivalent,
beirao2014hitchhiker, beirao2016virtual, de2016nonconforming, cangiani2017conforming}.
VEM is similar to the mimetic finite difference method (MFD) \cite{brezzi2005family, bochev2006principles, brezzi2009mimetic, cangiani2009convergence,
da2011arbitrary} in the sense of flexibility of mesh handling and using degrees of freedom only to construct the bilinear form.
However, MFD does not use basis functions while VEM assumes basis functions as solutions of local partial differential equations.
The word {\it virtual} comes from the fact that no explicit knowledge of the shape function is necessary.
By designing suitable elliptic projection operators on the local approximation space, VEM can be implemented using only the degrees of the freedom and the polynomial part of the approximation space, while the integration of source-term multiplied by virtual element test function on the right hand sides  is carefully handled using certain $L^2$-projection (see \cite{ahmad2013equivalent}).

The detailed guidelines for the implementation of VEM for elliptic problems including the construction of the projection operators can be found in
\cite{beirao2014hitchhiker, sutton2017virtual}.
Also,  nonconforming versions of VEM were studied in \cite{de2016nonconforming, cangiani2017conforming}.
 The developments and theories of VEM for elasticity problems and Stokes problems can be found in \cite{da2013virtual,gain2014virtual,da2015virtual,artioli2017stress,artioli2018family,zhang2019nonconforming} and  \cite{antonietti2014stream, cangiani2016nonconforming}, respectively.

On the other hand, the idea  of  VEM  was extended to the $H$(div) - conforming space on general polygons/polyhedral, called the mixed virtual element method (MVEM) in \cite{brezzi2014basic,da2016h,da2016mixed}, where the approximation spaces for the vector variables have degrees of freedom similar to those of BDM \cite{brezzi1985two} or Raviart-Thomas (RT) space \cite{RT1977}.

 The inner product term in the MVEM is defined  through an $L^2$-projection, thus the computations of the local integral is possible from the knowledge of degrees of freedom of elements, plus a stabilizing term which makes it compatible with ordinary inner product.
The MVEM leads to a saddle point problem similar to that of the  mixed finite element methods, which is a disadvantage of the mixed FEM.
Thus, it is necessary to devise a fast solution method for the algebraic equations arising from the mixed formulation of VEM. For example, an Uzawa type of solver may be used, or
a hybridization technique as in \cite{AB, brezzi1985two, dassi2021hybridization} can be employed. 
 Still, the resulting system  involving the Lagrange multipliers are nontrivial to solve; one has to invert the local matrix to find the  Schur complement.

In this paper, we develop new mixed VEM formulations for two and three dimensional problems along the line of mixed finite volume method (MFVM) introduced in \cite{Kwak2003,
kwak2012new}, where for the momentum equation, the gradient of test functions of a nonconforming space and some
subspace of polynomials are applied on each element, while the mass
equation is tested by a space of polynomials.
One of the advantages of the MFVM proposed in \cite{Kwak2003, kwak2012new} is that the formulation can be converted to the nonconforming finite element
method for the primary variable with modified forcing term.
Once the primary variable is obtained from solving the symmetric positive definite system, the velocity variable can be recovered locally.
Another advantage of this scheme is that the conservation of the momentum as well as the mass hold.

We develop a similar mixed  volume formulation using virtual elements on general polygonal/polyhedral meshes, by modifying the weak formulation introduced in \cite{kwak2012new}.
The $H$(div)-conforming VEM space in \cite{da2016h} or \cite{da2016mixed} is used for the vector variable, and the nonconforming VEM (NCVEM) space developed in \cite{de2016nonconforming} is used for the primary variable.

Our method is more naturally related to the NCVEM than MFEM is to nonconforming FEM, in the sense that the treatment of the forcing term is exactly the same as NCVEM
(i.e., one uses the $L^2$ projection on the right hand side.)

As is usual in VEMs, the variation form involves elliptic projection operators and stability terms for the primary variables, see (\ref{Newbil.CNC}), (\ref{box_scheme0}a).
By eliminating the velocity field from the first equation, we  obtain an equation for the NCVEM in the primary variable.
Once the primary variable is obtained by solving the (SPD) NCVEM system,
all the moments of the velocity variable can recovered locally.
Also, one can compute the $L^2$-projection of velocity variable easily.
Thus, the whole process can be implemented efficiently, avoiding the saddle point problems.
We name our method a {\bf mixed virtual volume method} (MVVM).

The proposed method is the first success in MVEM to compute the $H$(div) - conforming velocity variables by solving SPD problems in the
primary variable.
 Optimal error estimates for the proposed schemes are provided  and numerical results supporting our analysis are presented.
One may raise questions regarding the relationship of the proposed scheme with the reconstruction of velocity variable as in \cite{Marini}. Actually, the possibility is discussed in Section 4.2. In
the lowest order case,  we propose a one way to reconstruct an approximate velocity element of Raviart - Thomas type similar to \cite{Marini} in general
polygonal/polyhedral  mesh.

The rest of our paper is organized as follows.
The governing equation and brief review of  MVEM  are given in Section 2.
In Section 3, we review the nonconforming virtual element methods for the variable coefficient.
In Section 4, we introduce an MVVM and show that it is equivalent to the NCVEM.
 The error analysis is given in Section 5.
The numerical tests supporting our analysis are given in Section 6.
The conclusion follows in Section 7.

\section{Preliminaries}
 Let $\Omega$ be a bounded polygonal/polyhedral domain in $\R^d, d=2,3$ with the boundary $\bd\Omega$.
We consider the second-order elliptic boundary value problem
\begin{equation}\label{problem}
\left\{ \begin{aligned} -\div{\mathcal K}\Grad p   &= f \quad
\mbox{in } \Omega, \\
              p &= 0  \quad \mbox{on } \bd\Omega,
\end{aligned} \right.
\end{equation}
where $ {\cK}$ is a smooth, bounded, symmetric and uniformly positive definite tensor.

We introduce some notations  here: For any domain $D$, let $ H^{k}(D)$ (or $
\bH^{k}(D)$) be the scalar and vector Sobolev spaces of order $k\ge0$. We use
the standard notations $|\cdot|_{k,D}$,  $\norm{\cdot}_{k,D}$  and $\|\cdot\|_{\partial D}$  for
the (semi)-norms on $ H^{k}(D)$  and $L^2(\partial D)$, $(\cdot,\dot)_D$ for the $L^2$ inner product.  When $D=\Omega$, we drop the subscript
$\Omega$ and  write $|\cdot|_k, \norm{\cdot}_k$ instead.
 In two dimensions, we let
 $$\rot \bv=(\pd {v_2}x-\pd{v_1}y) \mbox{ and } \vrot q=(\pd {q}y,-\pd{q}x),$$
for  smooth enough vector and scalar functions $\bv$ and $q$.
 Let
\begin{eqnarray*}
 \bH(\div;D) &=& \{\bu \in (L^2(D))^d,  \mbox{$(d=2,3)$ with } \div\bu \in L^2(D)\},\label{Hdiv-local} \\
 \bH(\rot;D) &=& \{\bu \in (L^2(D))^2, \mbox{ with } \rot\bu \in L^2(D)\}, \label{Hcurl-local2D}\\
 \bH(\vcurl;D)&=& \{\bu \in (L^2(D))^3, \mbox{ with } \vcurl\bu \in (L^2(D))^3\}. \label{Hcurl-local3D}
\end{eqnarray*}
 The constants $C$, $C_*$ and $C^*$ will be independent of mesh size $h$, not necessarily the same for each occurrence.

 Let us introduce the vector variable $\bu = -{\mathcal K}\Grad p
$ and rewrite  problem $(\ref{problem})$ in the mixed form
\begin{equation}\label{mixed_system}
\left\{ \begin{aligned} \bu =- {\mathcal K}\Grad p &\quad \mbox{
in } \Omega, \\ \div\bu = f   &\quad \mbox{ in } \Omega,
\\ p = 0       &\quad \mbox{ on } \bd\Omega.
\end{aligned} \right.
\end{equation}
 Throughout this paper, we assume the following
regularity hold: The solution $(\bu,p)$ of (\ref{mixed_system})
satisfies $\bu\in \bH^{k+1}(\Omega)$, $p\in H^{k+2}(\Omega)$, and there
exists some constant $C>0$ such that
\begin{equation}\label{regul_system} \|\bu\|_{k+1} +\|p\|_{k+2}\le C\|f\|_k. \end{equation}
Its weak form is: find $\bu \in \bH(\div;D)$ and $p\in L^2(\Omega)$ such that
\begin{eqnarray} \label{mixed-conti}
({\mathcal K}^{-1} \bu, \bv)-(p,\mbox{div }\bv)&=&0 ,\quad \bv \in \bH(\div;\Omega), \\
 (\mbox{div } \bu,q) &=& (f,q), \quad q\in L^2(\Omega).
\end{eqnarray}

\subsection{Mixed virtual element methods}
We briefly review the mixed virtual element methods(MVEM)  introduced in
\cite{brezzi2014basic},\cite{da2016h},\cite{da2016mixed}.
Let $\mathcal{T}_h$ be a decomposition of $\Omega$ into regular polygons/polyhedra, and let $\cE_h^o$ be the set of all interior edges(faces), $\cE_h^\partial$ be the set of boundary edges(faces), and $\cE_h=\cE_h^o\cup \cE_h^\partial $.
Following \cite{da2016mixed}, \cite{de2016nonconforming}, we mean  by "regular" that, there exists some $\rho>0$ such that
\begin{itemize}
  \item $h_f \geq \rho h_\mathcal{P}$ holds for every element $\mathcal{P}\in \mathcal{T}_h$ and for every edge(face) $f \subset \partial \mathcal{P}$,
  \item every element $\mathcal{P}$ is star-shaped with respect to all points of a sphere of radius $\geq \rho h_\mathcal{P}$,
  \item when $d=3$, every face $f\in \mathcal{E}_h$ is star-shaped with respect to all points of a sphere of radius $\geq \rho h_{f}$,
\end{itemize}
where $h_f$(resp. $h_\mathcal{P}$) is the diameter of edge(face) $f$(resp. $\mathcal{P}$). We denote the maximum
diameter of elements $\Poly\in \cT_h$ by $h$.

For any  integer $k\ge0$, we denote by $P_k(D)$ the set of all polynomials of total degree less than or equal to $k$, and set
$P_{-1}(D)=\{0\}$. Also, we let the scaled polynomials:
\begin{equation}\label{Monomial}
 M_k(D)= \left\{ \left(\frac{\bx-\bx_D}{h_D} \right)^{\balpha},  |\balpha| \le k \right\},
\end{equation}  where  $\balpha=(\alpha_1,\cdots,\alpha_d)$ ($d=2,3$) is the multi-index and $\bx_D$  is the center of mass.

Let
\begin{eqnarray*} \label{Gradk}
 \cG_k(\Poly)&:=& \nabla P_{k+1}(\Poly),\\
 \, \cG_k(\Poly)^\perp & :=& \mbox{ orthogonal complement of } \cG_k(\Poly) \mbox{ in } (P_{k}(\Poly))^d,\\
 \cR_k(\Poly)&:=& \vcurl (P_{k+1}(\Poly))^3  \mbox{ if $d=3$ and $\rot(P_{k+1}(\Poly))^2$ if } d=2.\label{curl-rot}
\end{eqnarray*}
If we let $\pi_{k,d}$ be the dimension of $P_k(\RR^d)$, then we see
\begin{equation} \label{dimGradk}
 \mbox{dim\,}\cG_k(\Poly)=\pi_{k+1,d}-1, \,\mbox{dim\,} \cG_k(\Poly)^\perp= d \pi_{k,d}-\pi_{k+1,d}+1.
\end{equation}

Given $\Poly \in \mathcal{T}_h$,  the local  $\bH(\mathrm{div})$-conforming virtual element space  is defined as follows:
\begin{align} \label{localMVEM3D}
 \bV_{h}^k(\Poly):= \{&\bv\in \bH(\mdiv;\Poly) \cap \bH(\vcurl;\Poly) \, : \,  \bv\cdot \bn|_f \in P_k(f), \, \forall  \hbox{edges (faces) }  f \subset \partial \Poly, \nonumber \\
 & \mdiv \bv \in  P_{k}(\Poly), \,\vcurl\bv \in \cR_{k-1}(\Poly) \},
\end{align} where in two dimensional case, the $'\vcurl'$  operator is replaced by the  $'\rot'$ operator and  the space
$\cR_{k-1}(\Poly)$ is replaced by $P_{k-1}(\Poly)$.

 The global space of order $k$ is the space $\bV^k_h$ defined as
\begin{equation} \label{glabalMVEM}
 \bV^k_{h} = \{ \bv\in \bH(\mdiv;\OM): \bv|_{\Poly} \in  \bV^k_{h}(\Poly), \forall \Poly\in \cT_h \}.
\end{equation}
The  degrees of freedom  for $\bV_h^k$  are
\begin{eqnarray}
\frac1{|f|}\int_{f} \bv \cdot \bn g_k \,ds ,& &\,\, \forall g_k \in M_k(f),\, \forall f\in \cE_h,\label{Vdof3.7}\\
\frac1{|\Poly|}\int_{\Poly} \bv\cdot\bg_{k-1}\,dx ,& &\, \,\forall \bg_{k-1} \in \mathcal{G}_{k-1}(\Poly),\, \forall \Poly\in \cT_h, \label{Vdof3.8}\\
\frac1{|\Poly|}\int_{\Poly} \bv\cdot\bg_{k}^\perp\,dx ,& &\,\,\forall\bg_{k}^\perp\in\cG_k(\Poly)^\perp, \,\forall \Poly\in \cT_h.\label{Vdof3.9}
\end{eqnarray}
Here, $|\cdot|$ for any geometrical object means its Lebesgue measure and $g_{k}$, $\bg_{k-1}$,  $\bg_{k}^\perp$ are taken from the scaled monomials.
Let $\bPsi_h(\Poly)=\cG_{k-1}(\Poly) \oplus \cG_{k}^\perp(\Poly)$.
The conditions (\ref{Vdof3.8}), (\ref{Vdof3.9}) can be replaced by a single condition.
\begin{eqnarray*}
\frac1{|\Poly|}\int_{\Poly} \bv\cdot\bg \,dx ,& &\, \,\forall \bg  \in \bPsi_h(\Poly),\, \forall \Poly\in \cT_h.
\end{eqnarray*}

The pressure space is $$W_h^k:=\{q\in L^2(\Omega),\ q|_\Poly\in P_k(\Poly)\}.$$
\begin{remark} Let  $k\ge1$.
   Replacing the condition $\mdiv \bv\in P_k(\Poly)$ by $\mdiv \bv\in P_{k-1}(\Poly)$ in (\ref{localMVEM3D}) and replacing $k-1$ by $k-2$ in (\ref{Vdof3.8}), we obtain a BDM like virtual element space defined in \cite{da2016h}.
   However, we get $\mathcal{O}(h^k)$ instead of $\mathcal{O}(h^{k+1})$ in $H(div)$-norm. See Remark \ref{remark4.1} in Section 4.
\end{remark}
\subsection{Interpolations and $L^2$-projections}
The  $L^2$-projection operators $\Pi^0_k: L^2(\Poly) \rightarrow P_k(\Poly)$ and $\bPi^0_k:(L^2(\Poly))^d\rightarrow (P_k(\Poly))^d$ are defined as follows: On each $\Poly$, we define
\begin{equation}
\label{VL2proj}
\begin{cases}
\int_{\Poly} (q-\Pi_k^0 q) q_k \,dx=0,& \forall q_k\in P_k(\Poly),\\
\int_{\Poly} (\bv-\bPi_k^0 \bv) \bq_k \,dx=0,& \forall \bq_k\in (P_k(\Poly))^d.
\end{cases}
\end{equation}

When no confusion arises, we use the same notations $\Pi_k^0$ and $\bPi_k^0 $ to denote the $L^2$-projections from some virtual element spaces of $ L^2(\Omega)$ or $ (L^2(\Omega))^d$, although the computations are sometimes nontrivial (see the definition of
nonconforming virtual spaces in the next section).

As is shown in \cite{da2016h}, we can compute the $L^2$-projection $\bPi_k^0 \bv$ for $\bv\in \bV_h^k$  from the degrees of freedom of $\bv$ and the following properties hold: 
\begin{equation*} \label{Gradkttt}
 \|q-\Pi_k^0 q\|_0\le Ch^{k+1}|q|_{k+1},\quad \|\bv-\bPi_k^0 \bv\|_0\le Ch^{k+1}|\bv|_{k+1}.
\end{equation*}

The local interpolation operator $\bPi_k^F : (H^1(\Poly))^d\to \bV_h^k(\Poly)$ is defined by
\begin{eqnarray}
\int_{f} (\bv - \bPi_k^F \bv) \cdot\bn g_k \,d\sigma =0,& &\,\,  \forall g_k \in
M_k(f), \label{Vdof3.71}\\
\int_{\Poly} (\bv - \bPi_k^F \bv)\cdot \bg \,dx=0,& &\,\, \forall \bg \in \bPsi_h(\Poly). \label{Vdof3.81}
\end{eqnarray}

 Define  bilinear forms (for vector variables)
\begin{equation}
 \ba_h^\Poly(\bu,\bv):= (\cK\bPi^0_k\bu,\bPi^0_k\bv)_\Poly+\bS^\Poly(\bu-\bPi^0_k\bu,\bv-\bPi^0_k\bv)
\end{equation}
and
\begin{equation}
   {\ba}_h (\bu,\bv) =\sum_\Poly  {\ba}_h^\Poly(\bu,\bv),
\end{equation}
 where $\bS^\Poly(\bu,\bv)$ is any bilinear form that scales with the inner product $(\cK\cdot,\cdot)_\Poly$.

For $k\ge0$, the MVEM is : Find $ (\tilde{\bu}_h , \tilde{p}_h)\in \bV_h^{k}\times W_h^k $ such that
\begin{subequations} \label{VMixed}
\begin{gather} \ba_h(\tilde\bu_h, \bv_h)-(\tilde p_h,\mdiv \bv_h) = 0,  \quad \forall\bv_h\in \bV_h^{k}, \\
(\mdiv \tilde{\bu}_h, q_h)=(f,q_h), \quad \forall q_h \in W_h^k.
\end{gather}
\end{subequations}
The following error estimates are given in  \cite{da2016mixed}.
\begin{theorem}  \label{thm:2.1}
Under the assumptions above, the problem (\ref{VMixed}a,b) has a unique solution  $(\tilde\bu_h, \tilde{p}_h)$ and the following error estimates
hold.
\begin{eqnarray*}
 \|p-\tilde{p}_h\|_0  &\le& Ch^{k+1} (\|\bu\|_{k+1}+\|p\|_{k+1}),  \\
 \|\bu-\tilde\bu_h\|_{0} &\le & Ch^{k+1}\|\bu\|_{k+1}, \\
  \|\bu-\bPi_k^F\bu \|_{0} &\le & Ch^{k+1}\|\bu\|_{k+1}, \\
 \|\mdiv(\bu-\tilde\bu_h)\|_{0} &\le& Ch^{k+1} |f|_{k+1} .
\end{eqnarray*}

\end{theorem}
\section{Nonconforming virtual element methods}
We briefly describe the NCVEM introduced in \cite{de2016nonconforming},\cite{cangiani2017conforming}.

We need a broken Sobolev space
\begin{align*}
 H^1(\mathcal{T}_h)&=\left\{ q\in L^2(\Omega)  : q|_\Poly \in H^1(\Poly), \quad \forall \Poly \in \mathcal{T}_h \right\},
\end{align*}
with a broken norm
\begin{align*}
 \norm{q}^2_{1,h} = \sum_{\Poly \in \mathcal{T}_h} \norm{q}^2_{1,\Poly}.
\end{align*}
For positive integers $r=k+1, (k\ge0)$, we let
\begin{align}
 H^{1, nc}(\mathcal{T}_h;r)&=\left\{q\in H^1(\cT_h): \int_f  [q]m d\sigma =0,\, m\in P_{r-1}(f), \, \forall  f\in \cE_h^o  \right\} .
\end{align}%
In order to utilize the  nonconforming virtual element space in the next section, we need to use an extended version of VEM  as in   \cite{cangiani2017conforming}. The reason is to compute $L^2$-projection onto the space $P_r$.

The local space for NCVEM on each $\Poly \in \mathcal{T}_h$ is defined as
\begin{equation} \label{NoncLocalExt}
 N_h^r(\Poly)=\left\{q\in W_h^r(\Poly): (q-\Pi^*_r q,m)_{\Poly} =0,\, \forall m \in P_{r-1}(\Poly)\cup P_{r}(\Poly)\right\} ,
\end{equation}
where $ W_h^r(\Poly)$ is an  auxiliary space defined by
\begin{equation*} \label{NoncLocal}
 W_h^r(\Poly)=\left\{q\in H^1(\Poly): \pd q{\bn}\in P_{r-1}(f),\forall f\subset \partial\Poly, \Delta q \in P_{r}(\Poly)\right\},
\end{equation*}
and  $\Pi^*_r$ is a certain projection onto $P_r$ that can be computed from the degrees of freedom.
For example, one can use the elliptic projection $\Pi_r^\nabla$  \cite{beirao2013basic,cangiani2017conforming}.

 The global nonconforming virtual element space $N_h^r$ is defined as
\begin{equation} \label{glob.NCspace}
 N_h^r= \left\{q\in H^{1,nc}(\cT_h;r): q|_\Poly\in N_h^r(\Poly), \forall \Poly,  \,\int_f q m d\sigma =0, \, \forall m\in P_{r-1}(f),  \forall f \subset \cE_h^\partial \right\} .
\end{equation}
The global d.o.f.s are given by  the followings:
\begin{equation}\label{dofNC1}
\begin{aligned}
 \mu_{f,\balpha} (q) &=\frac1{|f|}\int_{f} q m_{\balpha}\, d\sigma, \forall m_{\balpha} \in M_{r-1}(f),\, f \in \cE_h^o ,\\
 \mu_{\Poly,\balpha} (q) &=\frac1{|\Poly|}\int_{\Poly} q m_{\balpha}\, dx, \forall m_{\balpha} \in M_{r-2}(\Poly),\, \Poly \in\cT_h.
\end{aligned}
\end{equation}
We define the usual elliptic  bilinear forms (for scalar variables) $a^\Poly: H^1(\Poly)\times H^1(\Poly) \rightarrow \R $ and  $a: H^1(\Omega)\times H^1(\Omega) \rightarrow \R $ as:
\begin{align*}
&a^\Poly(p, q ) = \int_\Poly \cK \nabla p \cdot \nabla q dx, \quad  \forall p,q \in H^1(\Poly), \\
&a(p,q) = \sum_{\Poly \in \mathcal{T}_h} a^\Poly(p,q), \quad \forall p,q \in H^1(\Omega).
\end{align*}
 Now we define a discrete bilinear form $a_h^\Poly(\cdot,\cdot): N_h^r\times N_h^r \to \RR$:
\begin{equation} \label{Newbil.CNC}
 a_h^\Poly(p_h,q_h)=(\cK\bPi_{r-1}^0 \nabla p_h,\bPi_{r-1}^0 \nabla q_h )_\Poly+ S^\Poly((I-\Pi_{r}^0)p_h),(I-\Pi_{r}^0) q_h),
\end{equation}
 where $S^\Poly$ is any stabilizing term satisfying
\begin{equation*} \label{bilinear.NCcoer}
 C_* a_h^\Poly(q_h,q_h) \le S^\Poly(q_h, q_h) \le C^* a_h^\Poly(q_h, q_h) , \forall q_h \in ker (\Pi_r^0).
\end{equation*}
We let
\begin{equation*} \label{globForm-a}
 a_h(p_h,q_h)= \sum_\Poly a_h^\Poly(p_h,q_h), \forall p_h, q_h \in N_h^{r} .
\end{equation*}
Now the NCVEM  of order $r\ge1$ is defined as  in \cite{cangiani2017conforming}:
Find  $p_h\in N_h^{r}$ such that
 \begin{equation} \label{NCVMFEM}
 a_h(p_h,q_h)= (\Pi_{r-1}^0f,q_h),
\end{equation}
The following optimal error estimate  for (\ref{NCVMFEM}) is given in Theorems 6.2 and 6.3 \cite{cangiani2017conforming}.
\begin{theorem}  Let $p$ and $p_h$ be  the solutions  of
(\ref{problem}) and (\ref{NCVMFEM}). Assume $p\in H^{r+1}(\OM),\ f\in H^{r-1}(\OM).$
Then, there exists a constant $C>0$ independent of $h$ such that
\begin{eqnarray*}
 \|p-p_h\|_{0}  +h|p-p_h|_{1,h} &\le& Ch^{r+1} \|p\|_{r+1}.
\end{eqnarray*}
 \end{theorem}

 \begin{remark}\label{remark3.1}
 \begin{enumerate}
   \item
 If we use $\Pi_{\max(r-2,0)}^0f$ on the right hand side of (\ref{NCVMFEM}), we can still get $H^1$ error estimate like
\begin{eqnarray*}
 |p-p_h|_{1,h} &\le& Ch^{r} (\|p\|_{r+1}+ \|f\|_{r-1}),
\end{eqnarray*}  but we do not get optimal $L^2$-error estimate.
\item
As is well known in VEM community, there are two choices of bilinear forms. We used the more general form (\ref{Newbil.CNC}) which works for variable coefficient. For constant coefficient $\cK$, the form using elliptic projection
$$(\cK\nabla \Pi^\nabla p_h,\nabla \Pi^\nabla q_h)+ S^\Poly((I-\Pi^\nabla)p_h),(I-\Pi^\nabla) q_h) $$ defined  in \cite{beirao2013basic,de2016nonconforming} can be used.
 \end{enumerate}
\end{remark}

\section{Mixed Virtual Volume Methods} Let $k\ge0$.
    Assume that we have  some $H(div)$ virtual  element space $\bV_h^{k}$  and NCVEM space  $N_h^{k+1}$ (to be associated with $\bV_h^{k}$).

We assume, for the sake of simplicity, that all the  elements $\Poly$ have  $n$ edges(faces), but our argument works when
each element has different number of edges(faces).
   We note the following type of Euler's formula:
\begin{equation} \label{EulerKim}
   n\#  \cT_h=\sum_{\Poly\in \cT_h}\sum_{f \subset \partial \Poly} 1=2\sum_{f\in \cE_h^o} 1+\sum_{f\in
 \cE_h^\partial}1=2\# \cE_h^o+\# \cE_h^\partial.
\end{equation}
%
Now we introduce our mixed virtual volume method (MVVM) for all order $k\ge0$:
Find $(\bu_h,p_h)\in\bV^k_h\times N_h^{k+1}$   which
 satisfies on every element $\Poly\in  \cT_h$,
\begin{subequations} \label{box_scheme0}
\begin{gather}\int_\Poly \bu_h \cdot\Grad\chi + a_{h}^\Poly (p_h,\chi) = 0,  \, \forall\chi\in N_h^{k+1}(\Poly), \\
\int_\Poly  (\bu_h +\cK \bPi^0_k\nabla p_h)\cdot \bv\,dx = 0,\, \forall \bv \in \cG_{k}(\Poly)^\perp,\ (k\ge1)  \\
\int_\Poly \div\bu_h \phi\,dx  = \int_\Poly \Pi_{k}^0 f \phi\,dx, \, \forall\phi\in P_k(\Poly) .
\end{gather}
\end{subequations}
From (\ref{box_scheme0}c) we have
\begin{equation}\label{divu=Pf}
 \mdiv  \bu_h = \Pi_k^0 f.
\end{equation}
 We see from (\ref{Vdof3.7}),(\ref{Vdof3.8}), (\ref{Vdof3.9}) (using (\ref{dimGradk})), that
  the dimension of $\bV^k_h\times N_h^{k+1}$  is
\begin{equation}\label{dimMVVMSp}
 \pi_{k,d-1}\#\cE_{h}\,+(\pi_{k,d}-1+d\pi_{k,d}- \pi_{k+1,d}+1)\#\cT_h + \pi_{k,d-1}\# \cE_h^o + \pi_{k-1,d} \#\cT_h
\end{equation}
while the number of equations in (\ref{box_scheme0}) is
\begin{eqnarray}
& (dim\, N_h(\Poly)-1+ dim\cG_{k}^\perp(\Poly)+dim\, P_k(\Poly))\cdot\# \cT_h\nonumber\\
=  \,&[ \pi_{k,d-1}n + \pi_{k-1,d}-1 + d\pi_{k,d}-\pi_{k+1,d}+1+ \pi_{k,d}]\# \cT_h.\label{dimMVVMeq}
\end{eqnarray}
Using Euler' formula (\ref{EulerKim}), we see
$$  \# \cE_{h}\,+  \# \cE_{h}^o=  n\#\cT_h .$$ Hence we see (\ref{dimMVVMSp}) and (\ref{dimMVVMeq}) are equal, and hence (\ref{box_scheme0}) is a square system.
Integration by parts gives
\begin{eqnarray}\label{eq:localIbyP}
  - \int_\Poly\bu_h\cdot\nabla\chi \,dx = - \int_{\partial \Poly}\bu_h\cdot\bn\chi\,ds +  \int_\Poly \mdiv\bu_h \chi\,dx .
\end{eqnarray}
Summing over all $\Poly$, we have, by (\ref{box_scheme0}a) and (\ref{divu=Pf})
\begin{eqnarray} %
   a_h (p_h,\chi)&=& -\sum_P  \int_\Poly\bu_h\cdot\nabla\chi\,dx\nonumber\\
   &=& -\sum_P\int_{\partial\Poly}\bu_h\cdot\bn\chi\,ds  +\sum_P \int_\Poly \Pi_k^0 f\chi\,dx.\label{eq:4.9}
\end{eqnarray}
 Now assume $\chi\in N_h^{k+1}$. Since $\chi$
  has continuous moments up to degree $k$ across internal edges(faces) and  has vanishing moments on $\partial \Omega$,
 we obtain
\begin{eqnarray}\label{equi-noncon}
a_h (p_h,\chi) &=&  (\Pi_{k}^0f,\chi),\quad \chi\in N_h^{k+1}.
\end{eqnarray}
This is exactly NCVEM  of order  $k+1$.
Thus we have shown that our mixed virtual volume scheme is equivalent to the  NCVEM.

\begin{remark}\label{remark4.1}

\begin{enumerate}
 \item For $k \ge 1$, we can
  replace  the test space in (\ref{box_scheme0}c) by $P_{k-1}(\Poly)$  to  obtain a scheme that corresponds to BDM like MVEM\cite{da2016h}, for which we lose one order in $H(div)$-norm.
  \item
  We can allow each polygon to have different number of edges (faces). 
  Similarly, the first term of the right hand side of (\ref{dimMVVMeq}) has to be changed exactly the same way. Thus the system is still a square system.
  \end{enumerate}
\end{remark}

\subsection{Recovery of $\bu_h$ and $L^2$-projection} %
We see  from (\ref{eq:localIbyP}), (\ref{divu=Pf}) and (\ref{box_scheme0}a)  that for any $\chi\in N_h^{k+1}(\Poly)$
\begin{equation}\label{eq:recov1}
\begin{aligned}
 \int_{\partial \Poly}\bu_h\cdot\bn \chi\,ds&=(\Pi_k^0f,\chi)_\Poly-a_h^\Poly(p_h,\chi).
\end{aligned}
\end{equation}
Hence the moments of $\bu_h\cdot\bn\in P_k(f)$  can be obtained
by choosing the basis functions
$\chi\in N_h^{k+1}$ corresponding to the degrees of freedom.

 The interior moments can be obtained similarly. 
 Indeed, for $\bv\in \cG_{k-1}(\Poly)=\nabla P_{k}(\Poly)$, we have
 $\bv=\nabla q$, for some $q\in P_{k}(\Poly)$. Hence
    from (\ref{box_scheme0}c), we have
\begin{eqnarray}\label{eq:RecG_k}
\int_\Poly \bu_h \cdot \nabla q\,dx =
 \int_{\partial\Poly}\bu_h\cdot\bn q\,d\sigma -\int_\Poly\Pi_k^0 f q\,dx,
\end{eqnarray} which is computable from the moments of $\bu_h\cdot\bn$.
 Meanwhile for $\bv\in \cG_k(\Poly)^\perp$, we see from (\ref{box_scheme0}b)
\begin{eqnarray}
\int_\Poly\bu_h \cdot\bv\,dx &=&-\int_\Poly\cK\bPi_{k}^0\nabla p_h\cdot\bv\,dx\nonumber \\
 &=&-\int_\Poly\nabla p_h\cdot \bPi_{k}^0(\cK\bv)\,dx \nonumber \\
&=&-\int_{\partial\Poly} p_h\bPi_{k}^0(\cK\bv)\cdot\bn\,d\sigma +\int_{\Poly} p_h \mdiv (\bPi_{k}^0(\cK\bv)) \,dx
\label{eq:RecG_kperp}
\end{eqnarray}
 which is computable from the d.o.f.s of $p_h$. Hence all the degrees of freedom of $\bu_h$ can be computed.

Furthermore, we can find the $L^2$-projection of $\bu_h$, at the same time.
 Since
\begin{equation*} \label{Grad_decomp2}
 (P_{k}(\Poly))^d=\cG_k(\Poly)\oplus  \cG_k(\Poly)^\perp,
\end{equation*}
and  $ \cG_k=\nabla P_{k+1}(\Poly)\subset \nabla N^{k+1}(\Poly)$, the same integration by parts  as
(\ref{eq:RecG_k}), (\ref{eq:RecG_kperp}) give  enough information to compute the projection
of $\bu_h$ onto $(P_k(\Poly))^d$.

 \subsection{Construction of Raviart - Thomas type approximation for $k=0$}
 It is clear that when $k=0$, MVVM and MFEM are equivalent on triangular/tetrahedral grids.
 In the case of MFEM, it is known \cite{Marini} that,  if $\cK$ is piecewise constant, then
\begin{equation} \label{Marini-inexpen} \bu_h= -\cK \nabla p_h +\frac{\bar f}d (\bx-\bx_\Poly), \quad d,=2,3,
\end{equation}
where $\bar f$ is the average of $f$ on each element  $\Poly$. The same formula holds for MVVM.

 On general grids, we have a similar representation (even when $\cK$ is nonconstant). We  see from (\ref{box_scheme0}c), that $\mdiv(\bu_h-\frac{\bar f}d (\bx-\bx_\Poly))=0$, $d=2,3$.
  Hence
   $$ \bu_h= \frac{\bar f}d (\bx-\bx_\Poly)+ \vcurl \xi(\vrot \xi \mbox{ if } d=2),\ \hbox{for some } \xi \in H^1(\Poly). $$
   Substituting this into   (\ref{box_scheme0}a), and letting $\chi=x, y,z$,  we have
    $$ \bPi_0^0(\vcurl \xi) = -\bPi_0^0(\cK \bPi_0^0( \nabla p_h))= -\bar \cK \bPi_0^0( \nabla p_h).$$
   Hence
  \begin{equation} \label{Marini-in2}
   \bu_h = \frac{\bar f}d (\bx-\bx_\Poly)-\bar \cK \bPi_0^0(\nabla p_h)+\vcurl\xi-\bPi_0^0(\vcurl\xi). \end{equation}
    The projection
    $ \bPi_0^0(\nabla p_h)$ can be computed by letting $\chi=x$ (or $\chi=y,z$) as follows:
 \begin{eqnarray*}
 \int_\Poly\bPi_0^0(\nabla p_h)\cdot\nabla \chi\,dx  &=&\int_\Poly\nabla p_h\cdot(1,0,0)^T\,dx = \int_{\partial\Poly} p_h n_x \,d\sigma\\
   &=& \sum_{f\subset\partial \Poly} n_x^f |f| \frac1 {|f|}\int_{f} p_h \,d\sigma =\sum_{f\subset\partial \Poly} n_x^f|f| \mu_{f,0}(p_h),
\end{eqnarray*}
 where $\mu_{f,0}(p_h)$ is the degree of freedom of $p_h$ on  $f$,
  $\bn=(n_x,n_y,n_z)$ is the unit outer normal to $\partial \Poly$, and
 $\bn^f=(n_x^f,n_y^f, n_z^f)$ is its restriction to each $f$.
 Hence
$$ \bPi_0^0 (\nabla p_h)=\frac1{|\Poly|}\begin{pmatrix} \sum_{f} n_x^f |f| \mu_{f,0} (p_h), & \sum_{f} n_y^f |f| \mu_{f,0} (p_h), &
\sum_{f} n_z^f |f| \mu_{f,0}(p_h)
\end{pmatrix}^T.$$
 We can write (\ref{Marini-in2}) in the form  (when $\bu_h$ is smooth)
 $$\bu_h =\tilde \bu_h+\vcurl\xi-\Pi_0^0(\vcurl\xi)=\tilde \bu_h +O(h),$$
 where
\begin{eqnarray}
 \tilde\bu_h & =&- \bar \cK\bPi_0^0 (\nabla p_h)+\frac{\bar f}d (\bx-\bx_\Poly) \in \bV_h. \label{RTlikescheme}
\end{eqnarray}
Thus we have obtained a lowest order Raviart - Thomas approximation to $\bu_h$ on polygonal/polyhedral grids.
 We believe it is a  better approximation than the  $L^2$-projection $\bar \cK\bPi_0^0(\bu_h)$, because it satisfies
 $\mdiv \tilde\bu_h =\bar f,$ while $\mdiv (\bar \cK \bPi_0^0(\bu_h))=0$.
 Indeed, the numerical tests support this assertion (see Table 6.2).

\section{Error estimates}

We need some lemmas which can be found in the literature.
\begin{lemma}
(Inverse inequality for VEM \cite{brenner2017some,chen2018some})
There exists a constant $C>0$ such that
\begin{eqnarray}\label{inv_ineq}
\norm{\nabla q}_0 \leq Ch^{-1} \norm{q}_0, \forall q \in N_h^{k+1} .
\end{eqnarray}
\end{lemma}
\begin{lemma}
(Norm equivalence for VEM \cite{brenner2017some,cangiani2017conforming,chen2018some})
For any $q \in N_h^{k+1}$, there exists a constant $C>0$ such that
\begin{eqnarray}\label{norm_equiv}
\frac{1}{C} h_\mathcal{\Poly}^{d/2} \norm{\mathbf{\Xi} (q)}_{\ell^2} \leq \norm{q}_{0,\Poly} \leq  C h_\mathcal{\Poly} ^{d/2} \norm{\mathbf{\Xi}(q)}_{\ell^2},
\end{eqnarray}
where $\mathbf{\Xi}(q)$ is the vector representing the degrees of freedom of $q$.
\end{lemma}

We  need the following lemma which are standard for FEM \cite{AB}, but not for VEM since there is no reference element.
\begin{lemma}\label{hybrid-recov}
Let $\phi \in L^2(\Poly)$,  $\bphi \in (L^2(\Poly))^d$ and $\mu\in L^2(\partial \Poly)$. Then
  the function
  $\chi\in N_h^{k+1}(\Poly)$ determined by
\begin{subequations}\label{Ndof}
\begin{align}
\int_f \chi q\,d\sigma&= \int_f \mu q\,d\sigma,\ \forall  q\in M_k(f), \mbox{ for all } f\subset\partial\Poly\\
\int_\Poly \chi m \,dx&= \int_\Poly  \phi m \,dx,\  \forall m\in M_{k-1}(\Poly)
 \end{align}
\end{subequations} satisfies
 \begin{equation}\label{NChi0}
 \|\chi\|_{0,\Poly} \le C (\|\phi \|_{0,\Poly}+ h^{1/2} \|\mu\|_{\partial \Poly}). \end{equation}
Similarly, the function $\bv\in \bV_h^{k}(\Poly)$ determined by the degrees of freedom
\begin{subequations}\label{A.udof}
\begin{align}
\int_f \bv \cdot\bn g\,ds&= \int_f \mu g\,ds,\, \forall g \in M_k(f), \, \mbox{ for all edges(faces) of } \Poly,\\
\int_\Poly \bv \cdot \mathbf{g}\,dx&= \int_\Poly \bphi \cdot \bg\,dx,\, \forall \bg\in \bPsi_h(\Poly)
\end{align}
\end{subequations}
  satisfies
  \begin{equation}  \label{A.Sigma0}
\|\bv\|_{0,\Poly} \le C(\|\bphi \|_{0,\Poly}+h^{1/2}\|\mu\|_{\partial \Poly}).
  \end{equation}
  \end{lemma}
 \begin{proof}
 We only prove (\ref{A.Sigma0}), since the proof of (\ref{NChi0}) is similar.
  It is well known that the square of $L^2$-norm of a function $\bv\in \bV_h^k(\Poly)$ scales like
 $$ |\Poly| \sum_{i=1}^{n_k} |\bXi_i(\bv)|^2,$$
 where  $\bXi_i(\bv)$ is  the $i$-th d.o.f of $\bv$ and $n_k$ is the dim $\bV_h^k(\Poly)$. In other words, if
 $\{\bphi_j\}_{j=1}^{n_k}$ is the canonical basis functions such that
 $\bXi_i(\bphi_j)=\delta_{ij}$, then
 $$ \bv = \sum_{i=1}^{n_k} \bXi_i(\bv)\bphi_j.$$
It can be easily verified that the scaled monomials $g\in M_k(f)$ and $\bg\in (M_k(\Poly))^d$ satisfy
\begin{equation}\label{A.norm_ofMonomial}
 \|g\|^2_{f}=O(|f|),\  \|\bg\|^2_{0,\Poly}=O(|\Poly|) .
\end{equation}
Hence for any $\bv\in\bV_h^k$, we have
\begin{eqnarray*} 
 \|\bv\|^2_{0,\Poly} & \le & C |\Poly|\sum_{i=1}^{n_k} |\bXi_i(\bv)|^2 \\
  & = & C |\Poly|[\sum (\mbox{edge(face) d.o.f.s})^2+ \sum(\mbox{interior d.o.f.s})^2]  \\
   &\le& C |\Poly| \left[\sum_{f\subset\partial\Poly,g}\left(\frac1{|f|}\int_f \bv\cdot\bn g\ d\sigma\right)^2 + \sum_{\bg}\left(\frac1{|\Poly|}\int_\Poly \bv\cdot\bg\ dx\right)^2 \right]  \\
  &\le & C  \frac{|\Poly|}{|f|^2}\sum_{f\subset\partial\Poly,g}\|\mu \|_{f}^2 \|g\|_{f}^2 +C |\Poly|^{-1}\sum_{\bg} \|\bphi\|_{0,\Poly}^2 \|\bg\|_{0,\Poly}^2\,(\mbox{by (\ref{A.udof})})  \\
  &\le & C \frac{|\Poly|}{|f|}\sum_{f\subset\partial\Poly,g}\|\mu\|_{f}^2 + C\sum_{\bg} \|\bphi\|_{0,\Poly}^2\, (\mbox{by (\ref{A.norm_ofMonomial})})\\
  & \le& C(h\|\mu\|_{\partial \Poly}^2+ \|\bphi\|_{0,\Poly}^2).
 \end{eqnarray*}
 \end{proof}

\begin{theorem}\label{thm:main}
Let $(\bu_h,p_h)$ be the solution of the system
$(\ref{box_scheme0})$. Then there exists a constant $C$ independent
of $h$ such that
\begin{subequations}\label{err2}
\begin{align}
   \norm{\bu- \bu_h}_0 & \le C h^{k+1}(\|\bu\|_{k+1} +\|f\|_{k}),  \\
 \norm{\div(\bu- \bu_h)}_0 &\le  Ch^{k+1} |f|_{k+1},
\end{align}
\end{subequations}   provided that $\bu\in \bH^{k+1}(\Omega)$ and $f\in  H^{k+1}(\Omega)$.
\end{theorem}
\begin{proof}
We shall prove  (\ref{err2}b) first. We see  from (\ref{divu=Pf})  that
 \begin{eqnarray*}
    \|\div(\bu-\bu_h)\|_0 &\le  \|f-\Pi_k^0f\|_0 \le Ch^{k+1}|f|_{k+1} .
\end{eqnarray*}
Next we prove  (\ref{err2}a).
  By the triangle inequality
$$ \| \bu- \bu_h\|_0 \le  \| \bu-\bPi_k^F\bu\|_0+ \|\bPi_k^F\bu-\bu_h\|_0, $$
and the approximation property of $\bPi_k^F$  (Theorem \ref{thm:2.1}),
it suffices to estimate
$$  \| \bPi_k^F\bu-\bu_h\|_0. $$
For the sake of simplicity we assume $\cK=1$.  (Similar estimate holds as long as $\cK$ is sufficiently smooth.)
  Let $p_\pi$ be an arbitrary function in $P_{k+1}(\Poly)$.
Then, clearly we have $a_h^\Poly(p_\pi,\chi)=a^\Poly(p_\pi,\chi)$ for all $\chi \in N_h^{k+1}(\Poly)$.
From (\ref{box_scheme0}a), we see
\begin{eqnarray}
 \int_\Poly (\bu-\bu_h)\cdot\Grad\chi\,dx &=& - a^{\Poly}(p,\chi ) + a_{h}^\Poly  (p_h,\chi) \nonumber\\
&=& - a^{\Poly}(p-p_\pi,\chi)+ a_h^\Poly (p_h-p_\pi,\chi).  \label{temp1}
\end{eqnarray}

 Let $\chi \in H^1(\Poly)$ be the solution of
  \begin{eqnarray*} \Delta \chi&=& 0 \text{ in }  \Poly,\\
  \chi &=&(\bPi_k^F \bu-\bu_h)\cdot\bn \text{ on } \partial\Poly.
 \end{eqnarray*}
 Since
 $\chi$ is completely determined by the moments $\int_f \chi q_k\,d\sigma = \int_f (\bPi_k^F \bu-\bu_h)\cdot\bn q_k\,d\sigma, \,\forall q_k\in P_k(f), \forall f\subset  \partial\Poly$, we see $\chi \in N_h^{k+1}(\Poly)$.
  Hence by (\ref{NChi0}) we have
 \begin{equation} \label{chi_Poly}
\|\chi\|_{0,\Poly} \le C h^{1/2} \| (\bPi_k^F \bu-\bu_h)\cdot\bn\|_{\partial\Poly} .
 \end{equation}
 Now by the definition of $\bPi_k^F$,  (\ref{temp1}), and (\ref{box_scheme0})
 \begin{eqnarray*}
  \int_{\partial\Poly} \chi^2\,d\sigma  &=&  \int_{\partial\Poly} \chi (\bPi_k^F \bu-\bu_h)\cdot\bn \,d\sigma\\
      &=&  \int_{\partial\Poly} \chi ( \bu-\bu_h)\cdot\bn \,d\sigma\\
  &=& \int_\Poly (\bu-\bu_h)\cdot\Grad \chi\,dx + \int_\Poly \mdiv(\bu-\bu_h) \chi\,dx\\
  &=&-a^{\Poly}(p-p_\pi,\chi)+ a_h^\Poly(p_h-p_\pi,\chi)+((I-\Pi_k^0)f,\chi)_\Poly.
\end{eqnarray*}
 Now by the  approximation property of $p_\pi$ and $p_h$, the inverse inequality, and (\ref{chi_Poly}),
\begin{eqnarray*}
  \|(\bPi_k^F \bu-\bu_h)\cdot\bn\|_{\partial\Poly}^2
  &\le& C (\|p-p_\pi\|_{1,\Poly}+\|p_h-p_\pi\|_{1, \Poly})\|\chi\|_{1, \Poly} +C h^{k}\|f\|_{k, \Poly}\|\chi\|_{0, \Poly} \\
  &\le& C (2\|p-p_\pi\|_{1,\Poly}+\|p_h-p \|_{1,\Poly})\|\chi\|_{1,\Poly} +C h^{k}\|f\|_{k, \Poly}\|\chi\|_{0, \Poly} \\
   &\le& C h^{k+1} \|p\|_{k+2,\Poly} \|\chi\|_{1,\Poly} +C h^{k}\|f\|_{k, \Poly}\|\chi\|_{0, \Poly}  \\
   &\le& C h^{k} (\|p\|_{k+2,\Poly}+\|f\|_{k, \Poly})\| \chi\|_{0, \Poly}, \ (\hbox{by the inverse inequality}) \\
   &\le& C h^{k}\|f\|_{k,\Poly} \|\chi\|_{0, \Poly}, \ (\hbox{by the regularity assumption}) \\
   &\le& C h^{k+1/2}\|f\|_{k, \Poly}\|(\bPi_k^F \bu-\bu_h)\cdot\bn\|_{\partial\Poly} .
\end{eqnarray*}
Hence
\begin{eqnarray}\label{A.2Error-edge}
  \|(\bPi_k^F \bu-\bu_h)\cdot\bn\|_{\partial\Poly}  \le  C h^{k+1/2} \|f\|_{k, \Poly}.
\end{eqnarray}
On the other hand,
 from equations (\ref{mixed_system}), (\ref{box_scheme0}a,b) we see
 \begin{eqnarray*}\label{testspbPsia}
\int_\Poly(\bu   +\cK \nabla p ) \cdot \bv\,dx &=&0, \, \forall\bv \in (L^2(\Poly))^d,  \\
\int_\Poly(\bu_h +\cK\bPi_k^0\nabla p_h) \cdot \bv\,dx &=&0, \, \forall\bv \in \bPsi_h(\Poly).\label{testspbPsib}
 \end{eqnarray*}
Subtracting, we have (since $\cK$ is constant)
 \begin{eqnarray*} \label{A.2intMom-sigma}
\int_\Poly (\bu -\bu_h)\cdot \bv\,dx &=& -\int_\Poly(\cK\nabla p - \cK\bPi_k^0 \nabla p_h)\cdot \bv\,dx,\,\, \forall\bv \in \bPsi_h(\Poly) \\
&=& -\int_\Poly(\cK\nabla p - \cK\nabla p_h)\cdot \bv\,dx,\,\, \forall\bv \in \bPsi_h(\Poly) .
 \end{eqnarray*}
Let $\mu=(\bPi_k^F \bu-\bu_h)\cdot\bn$ and $\bphi= \cK\nabla p - \cK\nabla p_h$. Then $ \bsigma=\bPi_k^F\bu-\bu_h \in \bV_h^k(\Poly)$ is the solution of
\begin{subequations}\nonumber 
\begin{align}
\int_f \bsigma \cdot\bn q\,ds&= \int_f \mu q\,ds,\, \forall q \in M_k(f), \, \mbox{ for all edges(faces) of } \Poly,\\
\int_\Poly \bsigma \cdot \mathbf{g}\,dx&= \int_\Poly \bphi(\Poly) \cdot \bg\,dx,\, \forall \bg\in
\bPsi_h.
\end{align}
\end{subequations}
Then by (\ref{A.Sigma0}), (\ref{A.2Error-edge}),  and the approximation property of $p_h$, we have
\begin{eqnarray*}
\|\Pi_k^F\bu-\bu_h\|_{0,\Poly}& \le &\|\cK\|_\infty  \| \nabla p -\nabla p_h\|_{0,\Poly} +Ch^{1/2}\|(\bPi_k^F \bu-\bu_h)\cdot\bn\|_{\partial\Poly}  \\
&\le&  C h^{k+1}(\|p\|_{k+2,\Poly}+\|f\|_{k,\Poly}). \label{A.Error-inte}
\end{eqnarray*}
By the triangle inequality and approximation property of $\Pi_k^F\bu$, the proof is complete.

\end{proof}

\begin{figure}[!ht]
\centering
\includegraphics[width = 0.42\textwidth]{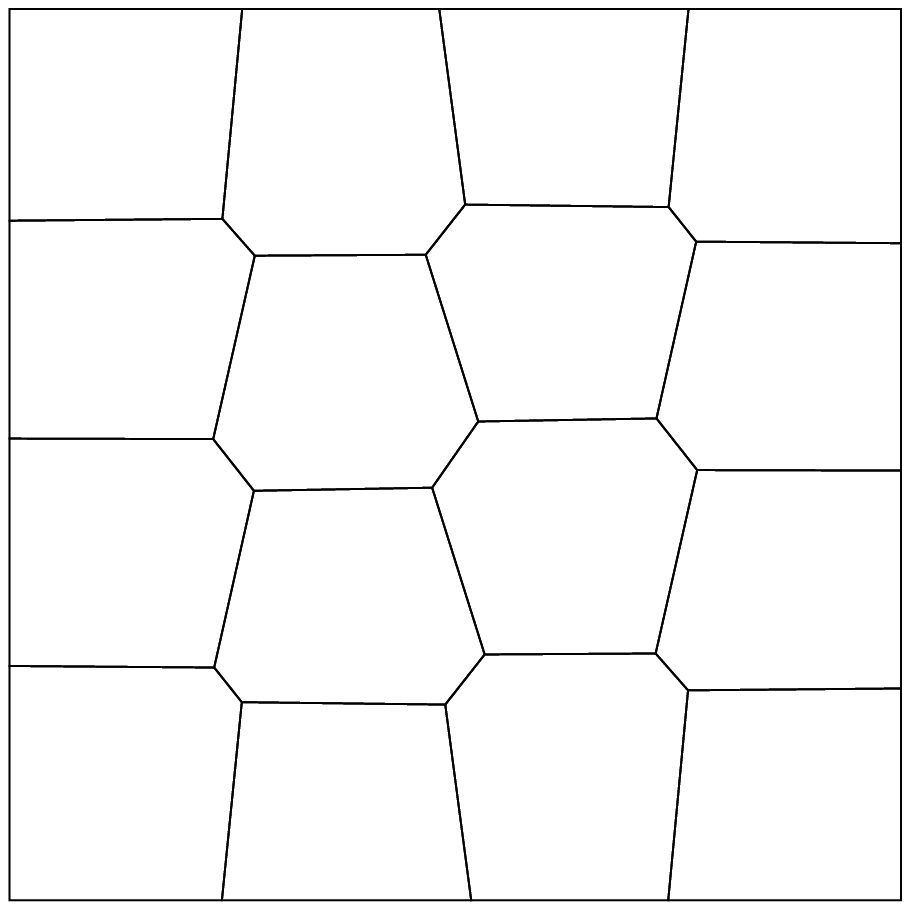}
\includegraphics[width = 0.42\textwidth]{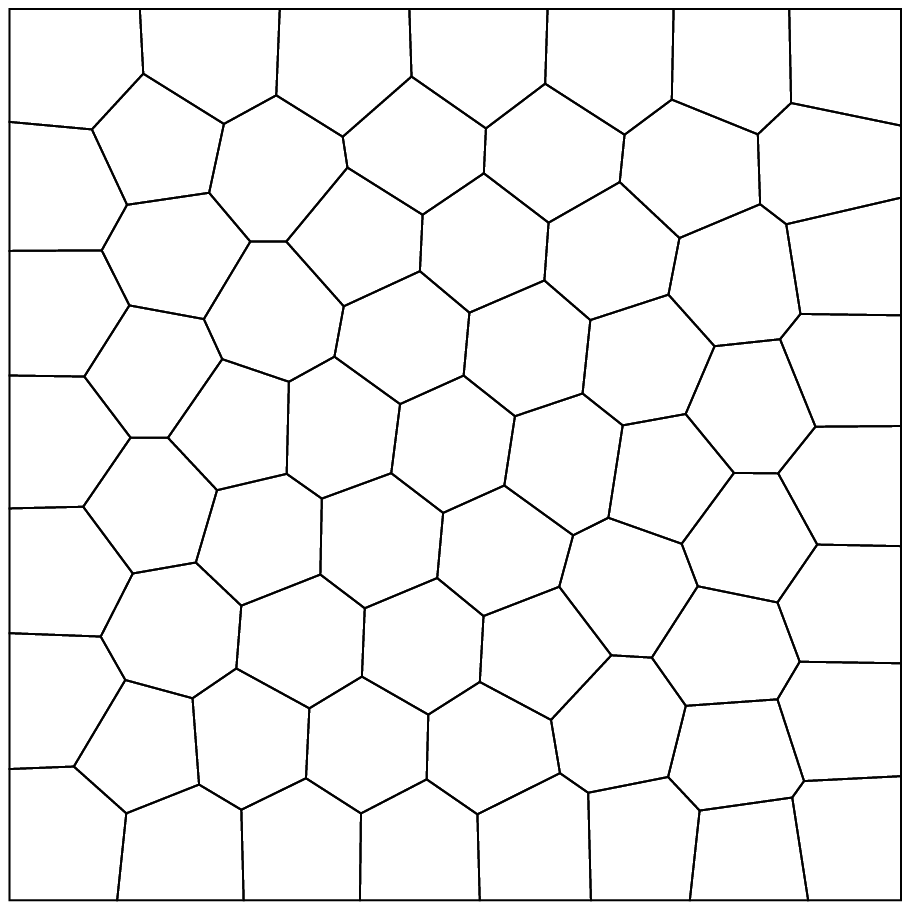}
\includegraphics[width = 0.42\textwidth]{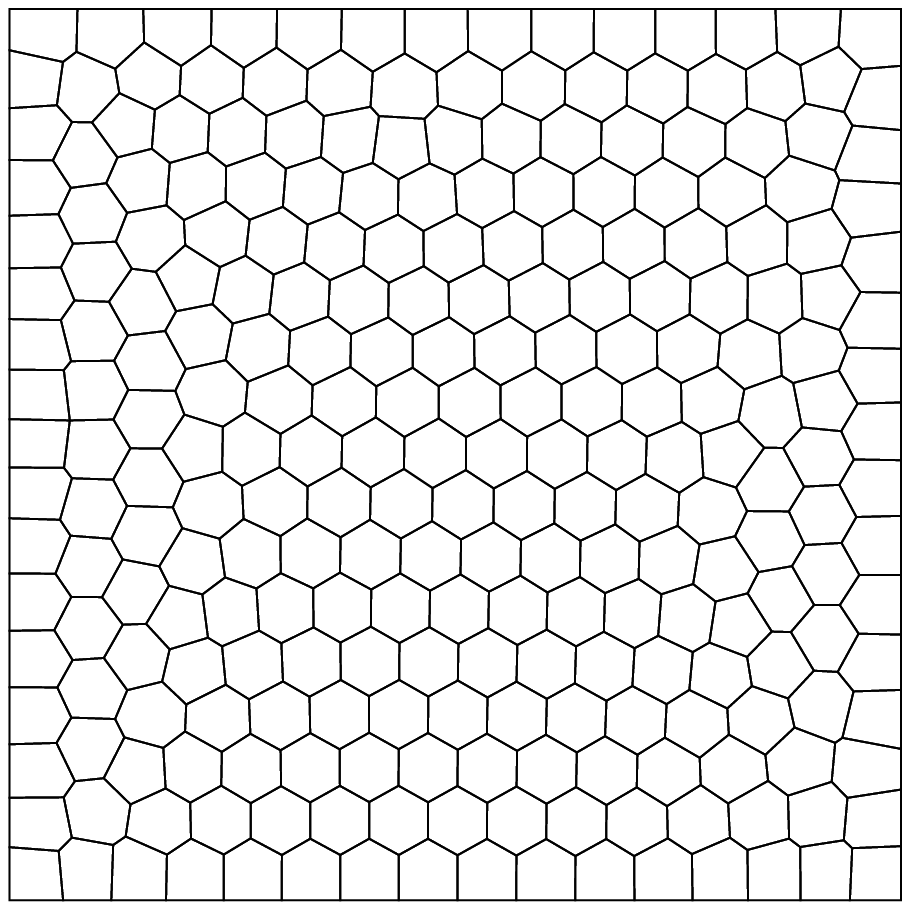}
\includegraphics[width = 0.42\textwidth]{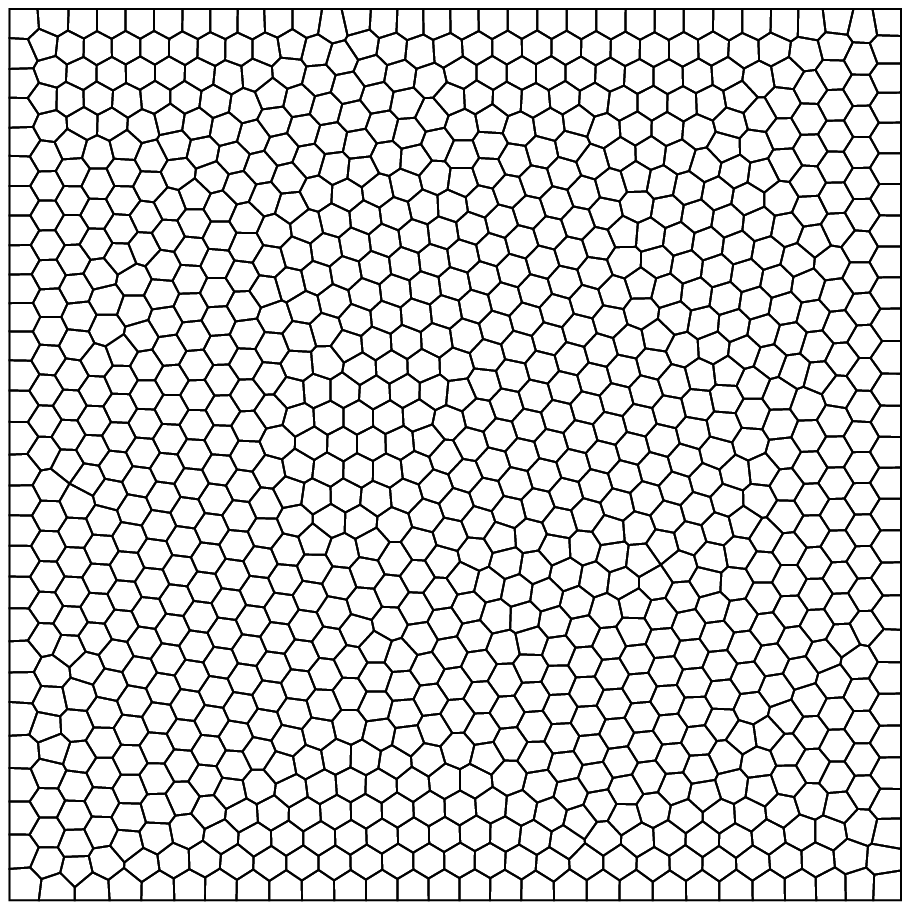}
\caption{Polygonal mesh $\mathcal{T}_h$ with number of elements 16 (left, top),  64 (right, top), 256 (left, bottom) and 1024 (right, bottom).}
\label{fig:mesh}
\end{figure}

\section{Numerical experiments}
In this section, we present some numerical results in two dimensional case.
The exact solutions on $\Omega=[0,1]^2$ are chose as
\begin{align*}
&p(x,y)   = x(1-x)y(1-y), \\
&\bu(x,y) =-\cK ( (1-2x)y(1-y), x(1-x)(1-2y)).
\end{align*}

  \begin{example} [Numerical results of MVVM]

We report the error between the exact solution and the $L^2$-projection of $\bu_h$ with $\cK = 1 + 0.5\sin(x)$, for $k=0,1,2,3$ in
  Table \ref{tab:allorder}.
We observe that the results are optimal for all cases.

\begin{table}[!ht]
\caption{$L^2$-errors between exact solution  and the $L^2$-projection of $\bu_h$ for orders $k=0,1,2,3$.
Table at left top, right top, left  bottom and right bottom correspond to the case $k=0$, $k=1$, $k=2$ and $k=3$ respectively .}
\centering
\begin{tabular}{r|rr}
\hline
$N$ of Elt. & $\|\bu- \bPi^0_0 \bu_h \|_0$  & order     \\
\hline
$2^2  \times 2^2$    &  6.494E-02         &       \\
$2^3  \times 2^3$    &  3.250E-02         & 0.998 \\
$2^4  \times 2^4$    &  1.578E-02         & 1.042 \\
$2^5  \times 2^5$    &  7.853E-03         & 1.007 \\
$2^6  \times 2^6$    &  3.895E-03         & 1.012 \\
\hline
\end{tabular}
\begin{tabular}{r|rr}
\hline
$N$ of Elt. & $\|\bu- \bPi^0_1 \bu_h \|_0$  & order \\
\hline
$2^2  \times 2^2$ & 2.964E-02         &       \\
$2^3  \times 2^3$ & 6.261E-03         & 1.981 \\
$2^4  \times 2^4$ & 1.162E-03         & 2.041 \\
$2^5  \times 2^5$ & 2.304E-04         & 2.001 \\
$2^6  \times 2^6$ & 5.040E-05         & 2.014 \\
\hline
\end{tabular}
\begin{tabular}{r|rr}
\hline
$N$ of Elt. &  $\|\bu- \bPi^0_2 \bu_h \|_0$  & order \\
\hline
$2^2  \times 2^2$  & 4.158E-03         &       \\
$2^3  \times 2^3$  & 6.260E-04         & 2.732 \\
$2^4  \times 2^4$  & 6.369E-05         & 3.297 \\
$2^5  \times 2^5$  & 5.719E-06         & 3.477 \\
$2^6  \times 2^6$  & 6.326E-07	       & 3.176 \\
\hline
\end{tabular}
\begin{tabular}{r|rr}
\hline
$N$ of Elt. &  $\|\bu- \bPi^0_3 \bu_h \|_0$  & order \\
\hline
$2^2  \times 2^2$  & 7.678E-05        &       \\
$2^3  \times 2^3$  & 8.684E-06        & 3.144 \\
$2^4  \times 2^4$  & 5.552E-07        & 3.967 \\
$2^5  \times 2^5$  & 3.368E-08        & 4.043  \\
$2^6  \times 2^6$  & 2.109E-09        & 3.997  \\
\hline
\end{tabular}
\label{tab:allorder}
\end{table}

 \end{example}
 \begin{example} [Comparison between $L^2$-projection and Raviart-Thomas type reconstruction]
In  the lowest order  case, we compare the errors of $L^2$-projection of $\bu_h$ and Raviart-Thomas type reconstruction (\ref{RTlikescheme}). 
Here, we set $\cK = 1$. The result   is reported in Table \ref{error_comparison}.
We observe that the Raviart-Thomas type is more accurate.
\begin{table}[!h]
\caption{$L^2$-errors of $\bPi^0_0 \bu_h$ (left) and Raviart-Thomas type reconstruction $\tilde\bu_h$ (right)}
\centering
  \begin{tabular}{c|c c|c c }
  \hline
   $N$ of Elt.  & $\|\bu - \bPi^0_0 \bu_h \|_{L^2(\Omega)}$   & order    & $\|\bu-\tilde \bu_h \|_{L^2(\Omega)}$  & order   \\
   \hline
    $2^2  \times 2^2$  & 4.303E-02	   &         & 3.171E-02       &          \\
    $2^3  \times 2^3$  & 2.241E-02     & 0.941   & 1.628E-02       & 0.962    \\
    $2^4  \times 2^4$  & 1.111E-02     & 1.012   & 7.930E-03       & 1.037   \\
    $2^5  \times 2^5$  & 5.575E-03     & 0.995   & 4.011E-03       & 0.983  \\
    $2^6  \times 2^6$  & 2.784E-03     & 1.002   & 1.988E-03       & 1.013 \\ \hline
  \end{tabular}
  \label{error_comparison}
\end{table}

 \end{example}

\section{Conclusion}
In this work, we  develop mixed virtual volume methods (MVVM) of all orders on polygonal/polyhedral meshes.
For the primary variable we use the nonconforming virtual element space, and for the velocity variable we use the $H(\mdiv)$ conforming   virtual element space.
The proposed method is the first success to compute $H$(div)-conforming velocity variables through the NCVEM.
We show that the  MVVM is equivalent to the NCVEM for all orders.
Once the primary variable is obtained from solving the (SPD) system arising from NCVEM, the velocity variable can be computed locally.
Thus, the whole procedure can be implemented efficiently, avoiding a saddle point problem.
The optimal error estimates are given and numerical results supporting our analysis are presented.


\end{document}